# Discussions On Driven Cavity Flow

Ercan Erturk

*Gebze Institute of Technology,
Energy Systems Engineering Department,
Gebze, Kocaeli 41400, Turkey*

SUMMARY

The widely studied benchmark problem, 2-D driven cavity flow problem is discussed in details in terms of physical and mathematical and also numerical aspects. A very brief literature survey on studies on the driven cavity flow is given. Based on the several numerical and experimental studies, the fact of the matter is, above moderate Reynolds numbers physically the flow in a driven cavity is not two-dimensional. However there exist numerical solutions for 2-D driven cavity flow at high Reynolds numbers.

Key Words: Driven Cavity Flow; Steady 2-D Incompressible N-S Equations; Fine Grid Solutions; High Reynolds Numbers

## 1  INTRODUCTION

The lid driven cavity flow is most probably one of the most studied fluid problems in computational fluid dynamics field. The simplicity of the geometry of the cavity flow makes the problem easy to code and apply boundary conditions and etc. Even though the problem looks simple in many ways, the flow in a cavity retains all the flow physics with counter rotating vortices appear at the corners of the cavity.

Driven cavity flow serve as a benchmark problem for numerical methods in terms of accuracy, numerical efficiency and etc. In the literature it is possible to find numerous studies on the driven cavity flow. The numerical studies on the subject of driven cavity flow can be basically grouped into three categories;

1) In the first category of studies, **steady** solution of the driven cavity is sought. In these type of studies the numerical solution of steady incompressible Navier-Stokes equations are presented at various Reynolds numbers.
2) In the second category of studies, the bifurcation of the flow in a driven cavity from a steady regime to an unsteady regime is studied. In these studies a **hydrodynamic stability analysis** is done and the Reynolds numbers at which a Hopf bifurcation occurs in the flow are presented.

---

Email: ercanerturk@gyte.edu.tr
URL: http://www.cavityflow.com

3) In the third category of studies, the transition from steady to unsteady flow is studied through a **Direct Numerical Simulation** (DNS) and the transition Reynolds number is presented.

In the first category studies that present steady solutions at high Reynolds numbers, Erturk et. al. [14], Erturk & Gokcol [15], Barragy & Carey [6], Schreiber & Keller [39], Benjamin & Denny [8], Liao & Zhu [30], Ghia et. al. [23] have presented solutions of *steady* 2-D incompressible flow in a driven cavity for Re≤10,000. Among these, Barragy & Carey [6] have also presented solutions for Re=12,500. Moreover, Erturk et. al. [14] and also Erturk & Gokcol [15] have presented steady solutions up to Re=20,000.

For the second category studies, Fortin et. al. [21], Gervais et. al. [22], Sahin & Owens [37] and Abouhamza & Pierre [1] are examples of *two dimensional* hydrodynamic stability studies on driven cavity flow.

In the third category, the following *two dimensional* Direct Numerical Simulation studies on driven cavity flow, Auteri, Parolini & Quartepelle [4], Peng, Shiau & Hwang [32], Tiesinga, Wubs & Veldman [44], Poliashenko & Aidun [34], Cazemier, Verstappen & Veldman [10], Goyon [24], Wan, Zhou & Wei [46] and Liffman [31] can be found in the literature as an example.

The above studies are all numerical studies. There are very few experimental studies in the literature on the driven cavity flow. Koseff & Street [27, 28, 29], Prasad & Koseff [35] have done several experiments on three dimensional driven cavity with various spanwise aspect ratios (SAR). These experimental studies present valuable information about the physics of the flow in a driven cavity.

Even though the driven cavity flow is studied at this extent in numerical studies, the nature of the flow at high Reynolds number is still not agreed upon. For example many studies from the first category, present steady solutions at very high Reynolds numbers showing that there exists a solution for steady 2-D Navier-Stokes equations for the flow inside a driven cavity. On the other hand, after a two dimensional hydrodynamic stability analysis or Direct Numerical Simulation, the studies from the second and third category claim that beyond a moderate Reynolds number the flow in a 2-D driven cavity is unsteady, therefore a steady solution **does not exist** hence a steady solution at high Reynolds numbers is not computable. These studies from second and third category basically contradict with the studies from the first category.

Based on the studies we found in the literature, we conclude that there are some confusions on the subject of driven cavity flow. We believe that for this flow problem some important points have to be discussed and analysed and agreed upon, such as;

- What is the physical nature of the flow in a driven cavity especially at high Reynolds numbers? Is the physical flow steady or turbulent or periodic in time? Is the flow two dimensional or three dimensional?
- What is the transition Reynolds number from a steady flow to unsteady flow?
- Is it possible to obtain a numerical solution for steady 2-D driven cavity at high Reynolds numbers? If it is possible then, does that mean that the real flow is steady and 2-D?
- What is the significance of a hydrodynamic stability analysis study of the flow in a driven cavity? In these studies, whether 2-D equations or 3-D equations should be used?

- What is the significance of a study that analyse the transition Reynolds number in a driven cavity through a Direct Numerical Simulation? In these studies, whether 2-D equations or 3-D equations should be used? How should we distinguish spurious solutions from real solutions?
- What is the significance of a study that presents 2-D steady solutions of driven cavity at high Reynolds numbers?

The main purpose of this study is then to discuss the incompressible flow in a 2-D driven cavity in terms of physical, mathematical and numerical aspects, together with a very brief literature survey on experimental and numerical studies. We will also present very fine grid steady solutions of the driven cavity flow at very high Reynolds numbers. Based on the results obtained, we will discuss on the driven cavity flow in an attempt to address the important points mentioned above.

## 2 STEADY SOLUTIONS OF 2-D DRIVEN CAVITY FLOW

In this section we will show that there exist numerical solutions for steady 2-D driven cavity flow at high Reynolds numbers and moreover these solutions are computable even with a simple explicit numerical method.

First, as the starting point, we assume that the incompressible flow in a driven cavity is **two-dimensional**. A second point is that, we are seeking a **steady** solution. After these two starting points, the problem at hand is governed by the 2-D steady incompressible Navier-Stokes equations. We use the governing equations in streamfunction ($\psi$) and vorticity ($\omega$) formulation such that

$$\frac{\partial^2 \psi}{\partial x^2} + \frac{\partial^2 \psi}{\partial y^2} = -\omega \tag{1}$$

$$\frac{1}{Re}\left(\frac{\partial^2 \omega}{\partial x^2} + \frac{\partial^2 \omega}{\partial y^2}\right) = \frac{\partial \psi}{\partial y}\frac{\partial \omega}{\partial x} - \frac{\partial \psi}{\partial x}\frac{\partial \omega}{\partial y} \tag{2}$$

where Re is the Reynolds number, and *x* and *y* are the Cartesian coordinates.

Erturk et. al. [14] have stated that for square driven cavity when fine grids are used, it is possible to obtain numerical solutions at high Reynolds numbers. At high Reynolds numbers thin boundary layers are developed along the solid walls and it becomes essential to use fine grid meshes. Also when fine grids are used, the cell Reynolds number or so called the Peclet number defined as $Re_c = u\Delta h/\nu$ decreases and this improves the numerical stability (see Weinan and Jian-Guo [47] and Tannehill et al. [43]). In this study, we use a very fine grid mesh with 1025×1025 grid points. On this mesh we solve the governing equations (1 and 2) using SOR method.

The vorticity values at the wall is calculated using Jensen's formula (see Fletcher [18])

$$\omega_0 = \frac{-4\psi_1 + 0.5\psi_2}{\Delta h^2} - \frac{3V}{\Delta h} \tag{3}$$

where subscript 0 refers to the points on the wall and 1 refers to the points adjacent to the wall and 2 refers to the second line of points adjacent to the wall, *V* refers to the velocity of the wall with being equal to 1 on the moving wall and 0 on the stationary walls and also $\Delta h$ is the grid spacing.

## 2.1 Results

As a measure of convergence to the steady state, during the iterations we monitored three residual parameters. The first residual parameter, RES1, is defined as the maximum absolute residual of the finite difference equations of steady streamfunction and vorticity equations (1 and 2). These are respectively given as

$$RES1_\psi = max\left(\left|\frac{\psi_{i-1,j}^{n+1} - 2\psi_{i,j}^{n+1} + \psi_{i+1,j}^{n+1}}{\Delta h^2}\right.\right.$$

$$\left.\left.\frac{\psi_{i,j-1}^{n+1} - 2\psi_{i,j}^{n+1} + \psi_{i,j+1}^{n+1}}{\Delta h^2} + \omega_{i,j}^{n+1}\right|\right)$$

$$RES1_\omega = max\left(\left|\frac{1}{Re}\frac{\omega_{i-1,j}^{n+1} - 2\omega_{i,j}^{n+1} + \omega_{i+1,j}^{n+1}}{\Delta h^2}\right.\right.$$

$$+\frac{1}{Re}\frac{\omega_{i,j-1}^{n+1} - 2\omega_{i,j}^{n+1} + \omega_{i,j+1}^{n+1}}{\Delta h^2}$$

$$-\frac{\psi_{i,j+1}^{n+1} - \psi_{i,j-1}^{n+1}}{2\Delta h}\frac{\omega_{i+1,j}^{n+1} - \omega_{i-1,j}^{n+1}}{2\Delta h}$$

$$\left.\left.+\frac{\psi_{i+1,j}^{n+1} - \psi_{i-1,j}^{n+1}}{2\Delta h}\frac{\omega_{i,j+1}^{n+1} - \omega_{i,j-1}^{n+1}}{2\Delta h}\right|\right) \quad (4)$$

The magnitude of RES1 is an indication of the degree to which the solution has converged to steady state. In the limit RES1 would be zero.

The second residual parameter, RES2, is defined as the maximum absolute difference between two iteration steps in the streamfunction and vorticity variables. These are respectively given as

$$RES2_\psi = max\left(\left|\psi_{i,j}^{n+1} - \psi_{i,j}^{n}\right|\right)$$
$$RES2_\omega = max\left(\left|\omega_{i,j}^{n+1} - \omega_{i,j}^{n}\right|\right) \quad (5)$$

RES2 gives an indication of the significant digit on which the code is iterating.

The third residual parameter, RES3, is similar to RES2, except that it is normalized by the representative value at the previous time step. This then provides an indication of the maximum percent change in $\psi$ and $\omega$ in each iteration step. RES3 is defined as

$$RES3_\psi = max\left(\left|\frac{\psi_{i,j}^{n+1} - \psi_{i,j}^{n}}{\psi_{i,j}^{n+1}}\right|\right)$$
$$RES3_\omega = max\left(\left|\frac{\omega_{i,j}^{n+1} - \omega_{i,j}^{n}}{\omega_{i,j}^{n+1}}\right|\right) \quad (6)$$

In our calculations, for all Reynolds numbers we considered that convergence was achieved when both $RES1_\psi \leq 10^{-10}$ and $RES1_\omega \leq 10^{-10}$ was achieved. Such a low value was chosen to ensure the accuracy of the solution. At these convergence levels the second residual parameters were in the order of $RES2_\psi \leq 10^{-17}$ and $RES2_\omega \leq 10^{-15}$, that means the streamfunction and vorticity variables are accurate to $16^{th}$ and $14^{th}$ digit accuracy respectively at a grid point and even more accurate at the rest of the grids.

Also at these convergence levels the third residual parameters were in the order of $RES3_\psi \leq 10^{-14}$ and $RES3_\omega \leq 10^{-13}$, that means the streamfunction and vorticity variables are changing with $10^{-12}$% and $10^{-11}$% of their values respectively in an iteration step at a grid point and even with less percentage at the rest of the grids. These very low residuals ensure that our solutions are indeed very accurate.

The boundary conditions and a schematics of the vortices generated in a driven cavity flow are shown in Figure 1 In this figure, the abbreviations BR, BL and TL refer to bottom right, bottom left and top left corners of the cavity, respectively. The number following these abbreviations refer to the vortices that appear in the flow, which are numbered according to size. Erturk et. al. [14] have presented an efficient numerical method and using their numerical method they have presented steady solutions of the cavity flow up to Reynolds number of 21,000 using a fine grid mesh. They have clearly stated that in order to obtain a steady solutions at high Reynolds numbers (Re>10,000), a grid mesh larger than 257×257 have to be used. In this study, first, we used a grid mesh with 257×257 grid points and solve the equations. With this many grid points we could not obtain a solution for Reynolds numbers above 10,000. Above $Re$=10,000 we observed that the solution was oscillating. Erturk et. al. [14] have stated that at high Reynolds numbers their solution was oscillating when they were using a coarse grid mesh. They also stated that when they increased the number of grids used, they were able to obtain a converged solution at the same Reynolds numbers. We, then, increased the number of grids to 513×513. This time we were able to obtain steady solutions up to Re≤15,000. With using this many number of grids, one thing is important to note, above $Re$=15,000 the solution was not converging but it was oscillating. Finally, we have increased the number of grids to 1025×1025. This time we were able to obtain converged steady solutions of driven cavity flow up to Reynolds number of 20,000. Above $Re$=20,000 our solution was oscillating again. This suggests that most probably, when a larger grid mesh is used the steady computations are possible.

Figure 2 and Figure 3 show the streamfunction and vorticity contours of the cavity flow up to Re≤20,000 with 1025×1025 grid mesh. These contour figures show that, the fine grid mesh provides very smooth solutions at high Reynolds numbers. The location of the primary and the secondary vortices and also the streamfunction ($\psi$) and vorticity ($\omega$) values at these locations are tabulated in Table 1. These results are in good agrement with Erturk et. al. [14]. They have stated that they have observed the quaternary vortex at the bottom left corner, BL3, appear in the solution at Re=10,000 when fine grids (600×600) were used. In this study we found that, when a more fine grid mesh (1025×1025) is used, BL3 vortex appear in the solution at Re=7,500. This would then suggest that in order to resolve the flow at high Reynolds numbers, fine grids are necessary.

Figure 4 and 5 show the u-velocity distribution along a vertical line and the v-velocity along a horizontal line passing through the center of the cavity respectively, at various Reynolds numbers. Detailed quantitative results, obtained using a fine grid mesh, was tabulated extensively in Erturk et. al. [14], therefore quantitative tabulated results will not be repeated in this study. The results obtained in this study agrees well with Erturk et. al. [14] and the reader is referred to that study for quantitative results.

We believe that, at high Reynolds numbers, numerical solutions needs to be validated since there is a chance that the solutions could be spurious. If there exists a solution to the 2-D steady incompressible equations at high Reynolds numbers, then this solution must satisfy the continuity of the fluid. The continuity will provide a very good mathematical check on the solution as it was suggested first by Aydin & Fenner [5]. We use the u- and v-velocity profiles in Figures 4 and 5 to test the accuracy of the solution. As seen in Figures 4 and 5 the integration of these velocity profiles will give the plus and minus areas shown by red colors. The degree to which the plus and minus areas cancel each other such that the integration give a value close to zero, will be an indicative of the mathematical accuracy of the solution. The velocity profiles are integrated using Simpson's rule to obtain the net volumetric flow rate $Q$ passing through these sections. The obtained volumetric flow rates are then normalized by a characteristic flow rate, $Q_c = 0.5$, which is the horizontal rate that would occur in the absence of the side walls (Plane Couette flow), to help quantify the errors. The obtained volumetric flow rate values ($Q_1 = |\int_0^1 u dy|/Q_c$ and $Q_2 = |\int_0^1 v dx|/Q_c$) are tabulated in Table 2. The volumetric flow rates in Table 2 are so small that they can be considered as $Q_1 \approx Q_2 \approx 0$. This mathematical check on the conservation of the continuity shows that our numerical solutions are indeed very accurate.

## 3  DISCUSSIONS ON DRIVEN CAVITY FLOW

The driven cavity flow problem has three aspects; physical (hydrodynamic), mathematical and numerical (computational) aspects and we think that the problem has to be discussed in terms of these aspects **separately**. We note that a very brief discussion on computational and also experimental studies on driven cavity flow can be found on Shankar & Deshpande [40]. First let us discuss the physical aspects of the cavity flow and look at the experimental studies.

Koseff & Street [27, 28, 29], Prasad & Koseff [35] have done several experiments on three dimensional driven cavity with various spanwise aspect ratios (SAR). Their experiments[27, 28, 29, 35] have shown that the flow in a cavity exhibit both local and global three-dimensional features. For example in a local sense they have found that Taylor-Görtler-Like (TGL) vortices form in the region of the Downstream Secondary Eddy (DSE). Also corner vortices form at the cavity end walls. In a global scale, due to no slip boundary condition on the end walls the flow is three-dimensional. They[27, 28, 29, 35] concluded that two-dimensional cavity flow **does not exist** (up to SAR=3) and further more, the presence of the TGL vortices **precludes** the possibility that the flow will be two dimensional even at large SAR. These conclusions clearly show that a two dimensional approximation for the flow in a cavity breaks down physically even at moderate Reynolds numbers.

Their results[27, 28, 29, 35] have also shown that the flow is laminar but unsteady even at moderate Reynolds numbers. The flow starts to show the signs of turbulence characteristics with turbulent bursts between 6,000≤Re≤8,000. Above Reynolds number of 8,000 the flow in a cavity is turbulent (for SAR=3).

According to these valuable physical information provided by [27, 28, 29, 35], the fact of the matter is, the flow in cavity is **neither two-dimensional nor steady** at high Reynolds numbers. Since two-dimensional steady flow in a cavity at high Reynolds numbers **does not** occur in reality, this flow (ie. 2-D steady cavity flow at high Re) is a **fictitious** flow, as it is also concluded by Shankar & Deshpande [40]. This is a very important fact to remember.

In our computations, as a starting point we have assumed that the flow inside a cavity is two-dimensional, therefore, we have used the 2-D Navier-Stokes equations and all of the solutions presented are based on the assumption that the flow is two-dimensional. This is also the case in most computational studies on driven cavity flow. At this point we should question whether or not an incompressible flow inside a cavity be two-dimensional at high Reynolds numbers. Let us consider a 3-D cavity with a moderate spanwise aspect ratio (SAR). In this cavity, the solution in z-axis will be reflection symmetric according to the symmetry plane (at the center in z-axis). For this given spanwise aspect ratio cavity, if the Reynolds number is small enough such that the effect of the no-slip conditions on the end walls in z-direction, on the flow at the symmetry plane is negligible, then the flow at the symmetry plane could be assumed two dimensional. However, in this cavity if the Reynolds number is high such that the end-wall-effects can no longer be assumed negligible, even though the flow is still reflection symmetric, the flow at the symmetry plane can not be assumed two dimensional since the velocity in the symmetry plane will not be 2-D divergence-free. In terms of the end-wall-effects only the flow in an infinite aspect ratio cavity (SAR $\rightarrow \infty$) could be assumed purely two dimensional at high Reynolds numbers and this is a fictitious scenario. We note that, here we look at the cavity flow only on a global scale and consider the end-wall-effects only. In local scale, experiments show that, TGL vortices appear in the region of the downstream secondary eddy (DSE) and [27, 28, 29, 35] state that the presence of the TGL vortices precludes the possibility that the flow will be two dimensional even at large SAR. At this point the study of Kim & Moin [26] is noteworthy to mention. They [26] have numerically simulated the three-dimensional time dependent flow in a square cavity with using periodic boundary conditions in the spanwise direction. They have observed TGL vortices in the flow field even though they did not have end walls in their simulations. Their results are important in the fact that in driven cavity flow the TGL vortices do not even need end walls to initiate. Therefore, even the flow in an infinite aspect ratio cavity will not be two-dimensional physically due to the TGL vortices. As a conclusion, physically, at high Reynolds numbers two dimensional cavity flow does not exist and any study that considers a two dimensional flow at high Reynolds numbers, is dealing with a fictitious flow.

Now let us discuss the cavity flow from a mathematical and numerical perspective. In our computations, we have used 2-D Navier stokes equations, and also we have presented **steady** solutions of the cavity flow at high Reynolds numbers. In light of these two points, what could be the nature of the 2-D incompressible flow in a driven cavity at high Reynolds numbers? Mathematically speaking a two-dimensional flow cannot be turbulent. Turbulence is by nature three-dimensional and is not steady. Therefore due to our two-dimensional assumption, mathematically the flow could not be turbulent. Since, when 2-D equations are used and therefore we do not let the flow be turbulent mathematically, could the flow still be time dependent, ie. periodic? Although 2-D high Reynolds number scenario is fictitious, this is a legitimate question in terms of a mathematical and numerical analysis. At this point we are trying to answer the mathematical and numerical nature of a fictitious flow and see if there exists a steady solution to 2-D Navier-Stokes equations or not. Mathematically speaking, a two-dimensional flow can be either steady (ie. solution is independent of time) or unsteady but time dependent (ie. solution is periodic in time) or unsteady but time independent (ie. solution is chaotic). Since there is a possibility that the flow could be either steady or unsteady, then is the flow in cavity steady or unsteady, or in other words does a steady solution to 2-D Navier-Stokes equations exist or not? In the literature it is possible to find studies that present steady solutions of 2-D incompressible Navier-Stokes equations at high Reynolds numbers and the studies of Erturk et. al. [14], Erturk and Gokcol [15], Barragy & Carey [6], Schreiber & Keller [39], Benjamin & Denny [8], Liao & Zhu [30], Ghia et. al. [23] can be counted as examples. In the literature there are also studies that claim the two dimensional cavity flow is unsteady, the studies in second and third category mentioned in the introduction can be counted as examples.

As mentioned earlier, in the second category studies, researchers have tried to obtain the Reynolds number at which a Hopf bifurcation occurs in the flow, ie. the Reynolds number at which the flow changes from steady to unsteady characteristics. In these studies, the basic steady cavity flow solution is perturbed with small disturbances and then the eigenvalues of linearized Navier-Stokes equations are analysed via hydrodynamic stability analysis.

Both external and internal flows, as Reynolds number is increased, exhibit a change from laminar to turbulent regime. In order to determine the Reynolds number that the transition from laminar to turbulent flow occurs, the Stability Theory is used. According to Stability Theory, in investigating the stability of laminar flows, the flow is decomposed into a basic flow whose stability is to be examined and a superimposed perturbation flow. The basic flow quantities are steady and the perturbation quantities vary in time. These basic and perturbation variables are inserted into Navier-Stokes equations. Assuming that the perturbations are small, equations are linearized and using a normal mode form for these perturbations, an Orr-Sommerfeld type of equation is obtained. The stability of a laminar flow becomes an eigenvalue problem of the ordinary differential perturbation equation. If the sign of imaginary part of the complex eigenvalue smaller than zero then the flow is stable, if it is greater than zero the flow is unstable (see Schlichting & Gersten [38] and also Drazin & Reid [13]).

In the case of driven cavity flow, the stability analysis requires solving the partial differential eigenvalue problem. Fortin et. al. [21], Gervais et. al. [22], Sahin & Owens [37] and Abouhamza & Pierre [1] are examples of hydrodynamic stability studies on driven cavity flow found in the literature. In these studies[21, 22, 37, 1] a **two-dimensional** basic flow is considered and then the solution of this two-dimensional basic flow is obtained numerically. Then this solution is perturbed with **two-dimensional** disturbances. The perturbation problem at hand is a partial differential eigenvalue problem, therefore the eigenvalues are also obtained numerically. These studies[21, 22, 37, 1] predict that a Hopf bifurcation takes place in a 2-D incompressible flow in a driven cavity some where around Reynolds number of 8,000.

One important point is that, the accuracy of the solution of the perturbation equation (or eigenvalues) completely depend on the accuracy of the solution of the basic flow. Since the perturbation quantities are assumed to be small compared to the basic flow quantities, any negligible numerical errors or oscillations in the basic flow solutions, will greatly affect the solution of the perturbation equation (Haddad & Corke [25], Erturk & Corke [16], Erturk et. al. [17]). Therefore for the sake of the accuracy of the solution of the perturbation equations, a highly accurate basic flow solution is required. In hydrodynamic analysis studies (second category studies) of driven cavity flow ([21, 22, 37, 1]) the number of grid points used is less than $257 \times 257$. Erturk et. al. [14] have reported that when a grid mesh with less than $257 \times 257$ points is used in cavity flow, the solution start to oscillate around Reynolds number range $7,500 \leq Re \leq 12,500$ depending on the order of the boundary conditions used at the wall. We note that the Reynolds numbers predicted by [21, 22, 37, 1] for Hopf bifurcation ($Re \approx 8,000$) is in the range, reported by Erturk et. al. [14] (7,500-12,500), where the solution becomes oscillatory due to coarse grid mesh (ie. large Peclet number). We believe that in these studies [21, 22, 37, 1] when the considered Reynolds number is close to the Reynolds number where the basic solution becomes oscillatory due to numerical instability, the solution of the perturbation equation is affected by this numerical instability in the basic flow. We also believe that in these studies if finer grids were used, the solutions would be different at high Reynolds numbers.

Most importantly, a hydrodynamic analysis would only be useful and also meaningful when there exist a physical flow otherwise the results of such a study would be physically useless. The experiments of Koseff & Street [27, 28, 29] and Prasad & Koseff [35] have shown that the flow inside a cavity is neither **two-dimensional** nor **steady** even at moderate Reynolds numbers. Therefore the results of a strictly **two-dimensional** hydrodynamic analysis (two-dimensional basic flow with two-dimensional perturbations) of a fictitious cavity flow will be also fictitious and will have no physical meaning.

Apart from the studies that examine the hydrodynamic stability of the driven cavity flow by considering a two-dimensional basic flow perturbed with two-dimensional disturbances mentioned above, the studies of Ramanan & Homsy [36] and Ding & Kawahara [12] are quite interesting and important. Both studies have considered a two-dimensional flow for driven cavity and analysed when this **two-dimensional** basic flow is perturbed with **three-dimensional** disturbances. Ramanan & Homsy [36] found that the two dimensional flow looses stability at Re=594 when perturbed with three-dimensional disturbances. Their results [36] seem to be in accordance with the experimental studies ([27, 28, 29, 35]). Also Ding & Kawahara [12] have predicted the instability at Re=1,025 with a frequency of 0.8. We believe that the results of Ding & Kawahara [12] over estimates the critical Reynolds number (Re=1,025) since they have used a slight compressibility in their simulations, such that their results are not purely incompressible. Note that the critical Reynolds numbers obtained with considering two-dimensional basic flow perturbed with three-dimensional disturbances in [12, 36] are much lower than the critical Reynolds numbers obtained with strictly two-dimensional hydrodynamic analysis in [21, 22, 37, 1] in which a two-dimensional basic flow perturbed with two-dimensional disturbances is considered. This is very important such that, the results of Ramanan & Homsy [36] and Ding & Kawahara [12] show that, for the hydrodynamic stability of driven cavity flow the spanwise modes are more dangerous than the two dimensional modes. Also note that, their results [36, 12] do not include the effect of the end walls in spanwise direction. We expect that with including the effect of end walls and the TGL vortices, a three-dimensional flow analysis will have a lower critical Reynolds number for stability when perturbed with three-dimensional disturbances.

As mentioned earlier, in the third category studies, researchers have tried to obtain the Reynolds number that the flow experience a transition from a steady regime to an unsteady regime using Direct Numerical Simulation (DNS). In these studies, first a Reynolds number is considered for computation and if a steady solution is obtained for this Reynolds number, then the flow is solved for a higher Reynolds number and this procedure is continued until a periodic solution is obtained. By doing several runs, the exact transition Reynolds number at which the solution changes characteristics from steady to unsteady behavior, is obtained. The following studies are example of DNS studies found in the literature.

Auteri, Parolini & Quartepelle [4] have used a second order spectral projection method. With this they have solved the Unsteady 2-D N-S equations in primitive variables. They have analysed the stability of the driven cavity with an impulsively started lid using 160×160 grids. They have removed the singularity that occurs at the corners of the cavity. They have increased the Reynolds number step by step until the solution becomes periodic. They have found that a Hopf bifurcation occurs in the interval Re=8,017.6-8,018.8. At this Reynolds number their solution was periodic with a frequency of 0.4496. They also reported that the cavity flow passes through a second Hopf bifurcation in the interval Re=9,687-9,765.

Peng, Shiau & Hwang [32] have used a Direct Numerical Simulation (DNS) by solving the 2-D unsteady Navier-Stokes equations in primitive variables. With using a maximum of 200×200 grids, they have solved the cavity flow by increasing the Reynolds number. At Re=7,402 ±4 their solutions became periodic with a certain frequency of 0.59. As Re was increased to 10,300 the flow became a quasi-periodic regime. When the Re was increased to 10,325 the flow returned to a periodic regime again. Between Reynolds numbers of 10,325 and 10,700 the flow experienced an inverse period doubling and between 10,600 and 10,900 the flow experienced a period doubling. Finally they claimed that the flow becomes chaotic when Re was greater than 11,000.

Tiesinga, Wubs & Veldman [44] have used Newton-Picard method and recursive projection method using a 128×128 grid mesh and with these they claimed that in a 2-D incompressible flow inside a cavity the first a Hopf bifurcation occurs at Re=8,375 with a frequency of 0.44. The next Hopf bifurcations occur at Re=8,600, 9,000, 9,100 and 10,000 with frequencies 0.44, 0.53, 0.60 and 0.70. At these intervals the flow is either stable periodic or unstable periodic.

Poliashenko & Aidun [34] have used a direct method based on time integration. Using Newton iterations with 57×57 number of grid points they have claimed that a Hopf bifurcation occurs in a lid driven cavity at Re=7,763 ±2% with a frequency of 2.86 ±1%. They have stated that this Hopf bifurcation is supercritical.

Cazemier, Verstappen & Veldman [10] have used Proper Orthogonal Decomposition (POD) and analysed the stability of the cavity flow. They compute the first 80 POD modes, which on average capture 95% of the fluctuating kinetic energy, from 700 snapshots that are taken from a Direct Numerical Simulation (DNS). They have stated that the first Hopf bifurcation take place at Re=7,819 with a frequency of about 3.85.

Goyon [24] have solved 2-D unsteady Navier-Stokes equations in streamfunction and vorticity variables using Incremental Unknowns with a maximum mesh size of 257×257. They have stated that a Hopf bifurcation appears at a critical Reynolds number between 7,500 to 10,000. They have presented periodic asymptotic solutions for Re=10,000 and 12,500.

Wan, Zhou & Wei [46] have solved the 2-D incompressible Navier-Stokes equations with using a Discrete Singular Convolution method on a 201×201 grid mesh. Also Liffman [31] have used a Collocation Spectral Solver and solved the 2-D incompressible Navier-Stokes with 64×64 collocation points. Both studies claimed that the flow in a cavity is periodic at Reynolds number of 10,000.

These are just example studies we picked from the literature that uses Direct Numerical Simulation (DNS) to capture the critical Reynolds number that a Hopf bifurcation occurs in a driven cavity flow.

For a moment, let us imagine that the 2-D incompressible flow inside a cavity is not stable at a given Reynolds number, such that a steady solution **does not** exist and the solution is **time dependent**, as predicted by [4, 32, 44, 34, 10, 24, 46, 31]. In this case, if one uses the **steady** N-S equations then he should not obtain a solution since there is no steady state solution. In previous section, we have shown that there exist steady state solutions for driven cavity at high Reynolds numbers and moreover, these steady state solutions are computable even with a simple explicit (SOR) method.

Now let us imagine the opposite, such that, the 2-D incompressible flow inside a cavity is stable at a given Reynolds number, thus the the flow is **steady** and a steady solution does exist. In this case, if one uses the **unsteady** N-S equations for computation, the solution should converge to this steady solution through iterations in time also. We note that, the boundary conditions are independent of time. Therefore, if there exists a unique solution to steady equations, then the solution of unsteady equations should converge to this steady solution also, since there is nothing to drive the solution vary in time because the boundary conditions are steady and also a steady solution of governing equations does exist. In this case, if the solution of unsteady equations do not converge to a steady solution, this would indicate that there are numerical stability issues with the solution of the unsteady equations and the time dependent periodic solutions or any other solution obtained that is different than the steady solution should not be trusted.

Erturk et. al. [14] have reported that at high Reynolds numbers when they used a grid mesh of 257×257, they observed that their solutions oscillate in the pseudo time. However, when they used a larger grid mesh than 257×257, they were able to obtain a steady solution at high Reynolds numbers. Similarly, in this study, when we used a grid mesh of 257×257, at high Reynolds numbers the solution was not converging to a steady state but it was oscillating, even when very small relaxation parameters were used. However when we increased the grid mesh to 513×513, we were able to obtain steady solutions up to Re=15,000. With this many number of grid points, above this Reynolds number, the solution started to oscillate again even when sufficiently small relaxation parameters that would not allow the solution to diverge were used. When we increased the grid points up to 1025×1025, again we were able to obtain steady solutions up to Re=20,000. Therefore based on the experiences of Erturk et. al. [14] and this study, we conclude that in order to obtain a steady solution for the driven cavity flow, a grid mesh larger than 257×257 is necessary when high Reynolds numbers are considered and also at high Reynolds numbers when a coarse grid mesh is used then the solution oscillates. The interesting thing is that, both in this study and in Erturk et. al. [14], the obtained false periodic numerical solutions looked so real and fascinating with certain frequencies and periodicity. We believe that the studies that presented unsteady solutions of driven cavity flow using Direct Numerical Simulations ([4, 32, 44, 34, 10, 24, 46, 31]) have experienced the same type of numerical oscillations because they have used a small grid mesh. We would like to note that in all of the Direct Numerical Simulation studies on the driven cavity flow found in the literature ([4, 32, 44, 34, 10, 24, 46, 31]), the maximum number of grid points used is 257×257. We believe that because of the course grids used in [4, 32, 44, 34, 10, 24, 46, 31], their periodic solutions resemble the false periodic solutions observed in Erturk et. al. [14] and also in this study when a 257×257 coarse grid mesh is used. We also believe that, if a sufficiently fine grid mesh is used, a Direct Numerical Simulation (DNS) algorithm would also confirm the same **steady** results obtained in this study and in Erturk et. al. [14] at high Reynolds numbers. As Erturk et. al. [14] have stated, one of the reasons why the steady solutions of the driven cavity flow at very high Reynolds numbers become computable when finer grids are used, may be the fact that as the number of grids used increases, ie. $\Delta h$ gets smaller, the cell Reynolds number or so called Peclet number defined as $Re_c = \frac{u \Delta h}{v}$ decreases. This improves the numerical stability characteristics of the numerical scheme used (see [43, 47]), and allows high cavity Reynolds number solutions computable. Another reason may be that fact that finer grids would resolve the corner vortices better. This would, then, help decrease any numerical oscillations that might occur at the corners of the cavity during iterations. When we used a coarse grid mesh in our computations, the oscillations we observed looked so real that they could easily be mistaken deceptively as the real time behavior of the flow field. One may think that if the numerical simulation of time dependent equations (DNS) do not converge to steady state, then the flow is not hydrodynamically stable. However one should not forget that a DNS is also restricted with certain numerical stability conditions, such as a Peclet number restriction in driven cavity flow as explained. We believe that, in a DNS study extensive grid study should be done to verify the results otherwise any numerical oscillations in the solution due to numerical stability issues can easily be confused as the hydrodynamic oscillations.

Among the Direct Numerical Simulation (DNS) studies found in the literature we believe that the studies of Albensoeder et. al. [2] and Albensoeder and Kuhlmann [3] are the most important studies on the driven cavity flow. In these studies, Albensoeder et. al. [2] and Albensoeder and Kuhlmann [3], have solved the three dimensional Navier Stokes equations using a Chebyshev-collocation method, therefore we believe that their [2, 3] solutions are very accurate. According to [2, 3] the steady two-dimensional flow in an infinite Spanwise Aspect Ratio (SAR) cavity with a square cross-section becomes unstable to steady short-wavelength Taylor–Görtler vortices at Re=786.3 ±6. We note that this 3-D DNS solution of critical bifurcation Reynolds number (Re=786) obtained by Albensoeder et. al. [2] and Albensoeder and Kuhlmann [3] is almost an order lower than the 2-D DNS solutions of critical bifurcation Reynolds number obtained by [4, 32, 44, 34, 10, 24, 46, 31].

The two dimensional cavity flow at high Reynolds number is a fictitious flow. Mathematically and numerically it is possible to study fictitious flows and in the literature it is possible to find many fictitious flows that are subject of mathematical and numerical studies. For example, the 2-D incompressible flow over a circular cylinder is steady for Reynolds number up to approximately 40. Beyond that Reynolds number, there appears Karman vortex street at downstream of the cylinder and the physical, ie. the real flow is unsteady. A steady solution beyond Re=40, if exists, is fictitious. However for this flow case, the 2-D incompressible flow over a cylinder, it is possible to obtain a **steady** solution mathematically when Re goes to infinity as the limiting solution (Kirchoff-Helmholtz solution, see Schlichting & Gersten [38]). Smith [41] and Peregrine [33] have done detailed mathematical analysis on 2-D **steady** incompressible flow over a circular cylinder at Reynolds numbers higher than 40. Apart from these mathematical studies on a fictitious flow, ie. 2-D steady incompressible flow over a circular cylinder above Re=40, in the literature there are also numerical studies on the same flow problem. Fornberg [19, 20], Son & Hanratty [42], Tuann & Olson [45] and Dennis & Chang [11] have presented numerical solutions of steady flow past a circular cylinder at Re=300, 600, 500, 100 and 100 respectively, which are much larger than Reynolds number of 40. This fact is important such that, even though the flow is physically fictitious, a mathematical solution exists and also it is possible to obtain a numerical solution as well.

In a mathematical study, Burggraf [9] applied Batchelor's model [7], which consists of an inviscid core with uniform vorticity coupled to boundary layer flows at the solid surface, to steady incompressible 2-D driven cavity flow at high Reynolds numbers. Burggraf [9] analytically calculated the theoretical core vorticity value at infinite Reynolds number as 1.886. As stated in Erturk et. al. [14], when the Reynolds number increases, thin boundary layers are developed along the solid walls and the core fluid moves as a solid body with a uniform vorticity in the manner suggested by Batchelor [7]. The numerical solutions of Erturk et. al. [14] agree well with the analytical solution of Burggraf [9] such that, as the Reynolds number increases the computed vorticity value at the primary vortex asymptotes to the theoretical infinite Re vorticity value.

Figure 6 compares the computed vorticity value at the center of the primary vortex with the theoretical infinite Reynolds number core vorticity value. In this figure the dotted line shows the theoretical vorticity value. In order to see the effect of the grid spacing on the vorticity value at the core of the primary eddy, in this figure looking at the solutions of Erturk et. al. [14] with 401×401, 513×513 and 601×601 grid mesh (tabulated in Table III, page 757 in [14]) and the solutions of this study with 1025×1025 grid mesh, it is clear that when a coarse grid mesh is used (for second order spatial accuracy, $\mathcal{O}(\Delta h^2)$) the computed vorticity value is less than the theoretical value in absolute value at high Reynolds numbers and as the number of grids increases the computed vorticity value gets closer to the theoretical value. The computed value should asymptote to the theoretical value while remaining greater than the theoretical value in absolute value. Therefore, for second order spatial accurate solutions, $\mathcal{O}(\Delta h^2)$, of driven cavity flow even 1025×1025 number of grids can be considered as a coarse grid mesh at very high Reynolds numbers. In order to see the effect of the spatial accuracy on the vorticity value at the core of the primary eddy, in this figure comparing second order $\mathcal{O}(\Delta h^2)$ 601×601 grid mesh solutions of Erturk et. al. [14] with fourth order $\mathcal{O}(\Delta h^4)$ 601×601 grid mesh solutions of Erturk & Gokcol [15], we see that as the spatial accuracy increases the computed value agree better with the theoretical value. This figure clearly shows that for 2-D steady incompressible driven cavity flow, at high Reynolds numbers higher order approximations together with the use of fine grids are necessary for accuracy.

## 4 CONCLUSIONS

In this study, the steady incompressible 2-D driven cavity flow is discussed in detail in terms of physical and mathematical and also numerical aspects. After a brief literature survey and a detailed discussion on the physical, mathematical and numerical characteristics of the incompressible flow in a driven cavity, we conclude the followings.

1. Physically, the flow in a driven cavity is neither two-dimensional nor steady, most probably, even at Re=1,000.
2. At high Reynolds numbers, when the incompressible driven cavity flow is considered as two-dimensional and also steady, then the considered flow is a fictitious flow.
3. It would be needless to study the hydrodynamic stability of a fictitious flow, ie. the 2-D steady incompressible flow in a driven cavity at high Reynolds numbers. Bifurcation Reynolds numbers obtained using a three dimensional hydrodynamic stability analysis differ an order from that obtained using a two dimensional hydrodynamic stability analysis of the driven cavity flow problem.
4. When Direct Numerical Simulation (DNS) is used in order to obtain the hydrodynamic stability (ie. the critical Reynolds number) of the driven cavity flow problem, extensive grid study should be done and great care should be given to distinguish any physical behavior from spurious behavior. 3-D DNS solutions for bifurcation Reynolds number differ an order from that obtained from 2-D DNS solutions of the driven cavity flow problem.
5. Mathematically, it is always possible to obtain a steady solution of a fictitious flow at the limiting case when Re goes to infinity. As shown in Erturk et. al. [14], Bachelor's [7] model of recirculating flow confined in closed streamlines at infinite Reynolds number seems to work with the square driven cavity flow. The numerical solutions of Erturk et. al. [14] agrees well with the analytical solutions of Burggraf [9].
6. Numerically, it is possible to obtain numerical solutions of 2-D steady incompressible cavity flow at high Reynolds numbers when fine grid meshes are used.
7. The model flow problem, the 2-D steady incompressible driven cavity flow, serve as a good benchmark problem for different numerical methods and boundary conditions, in terms of accuracy, convergence rate and etc., provided that these numerical solutions should be used for numerical comparison purposes between different solutions or with the analytical solution.
8. The fact that we can obtain numerical solutions of a particular 2-D steady flow problem does not necessarily mean that the actual physical flow is two-dimensional and steady. Similarly, the fact that a particular physical flow problem is not two-dimensional and steady does not necessarily mean that we cannot obtain 2-D steady numerical solutions of the particular flow.

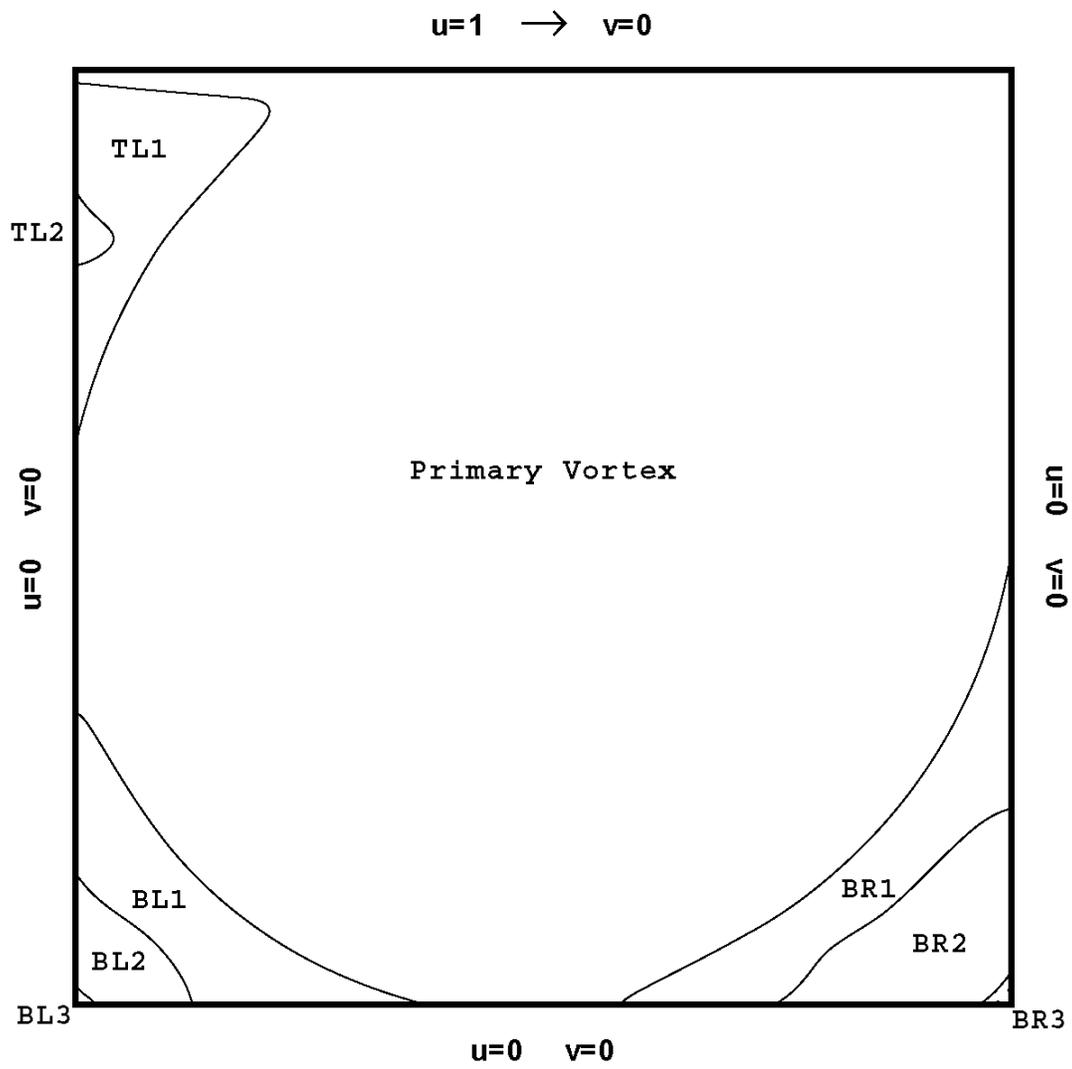

Figure 1. Schematic view of driven cavity flow.

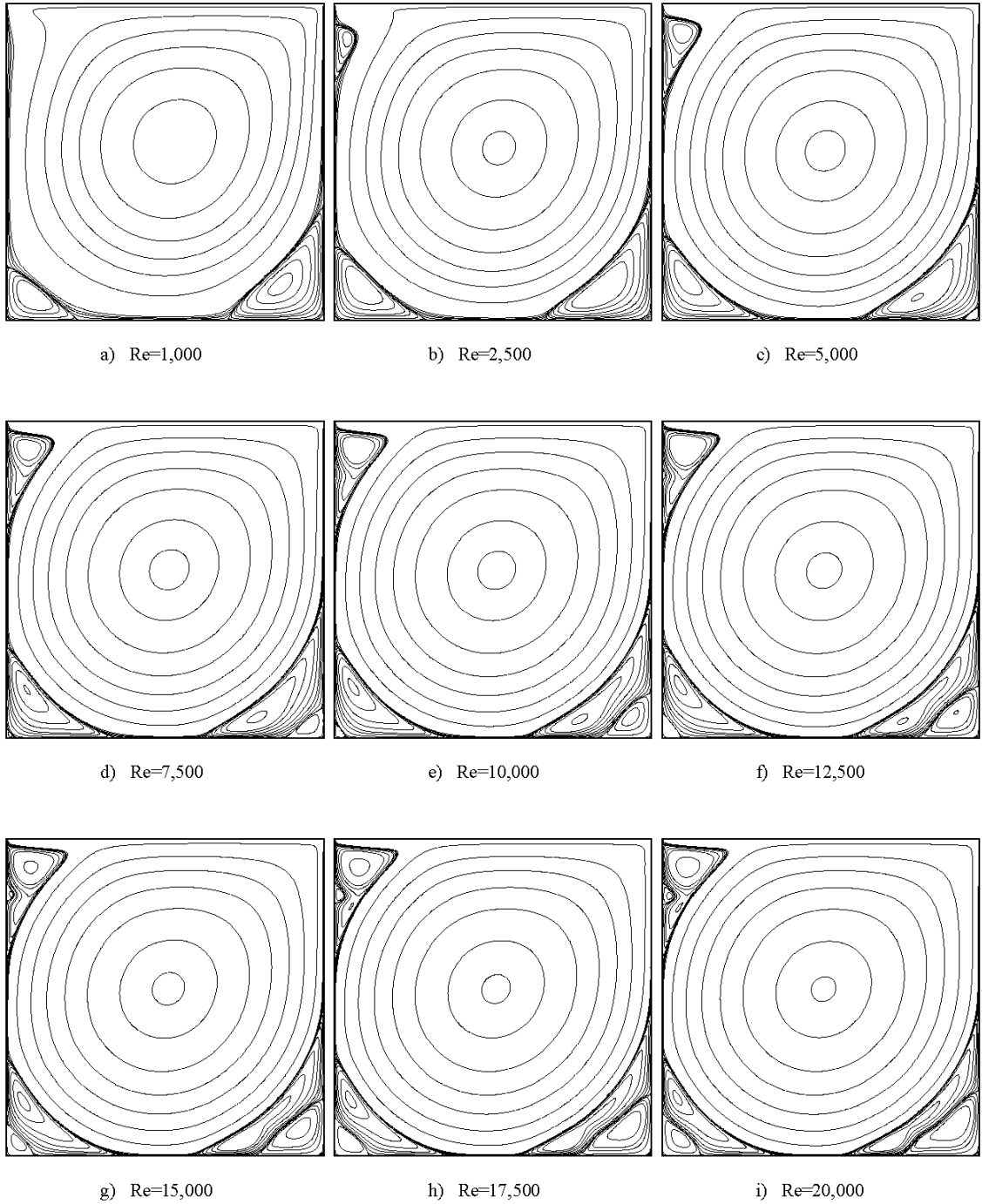

Figure 2. Streamfunction contours at various Reynolds numbers: (a) Re=1000;
(b) Re=2500; (c) Re=5000; (d) Re=7500; (e) Re=10000; (f) Re=12500;
(g) Re=15000; (h) Re=17500; and (i) Re=20000.

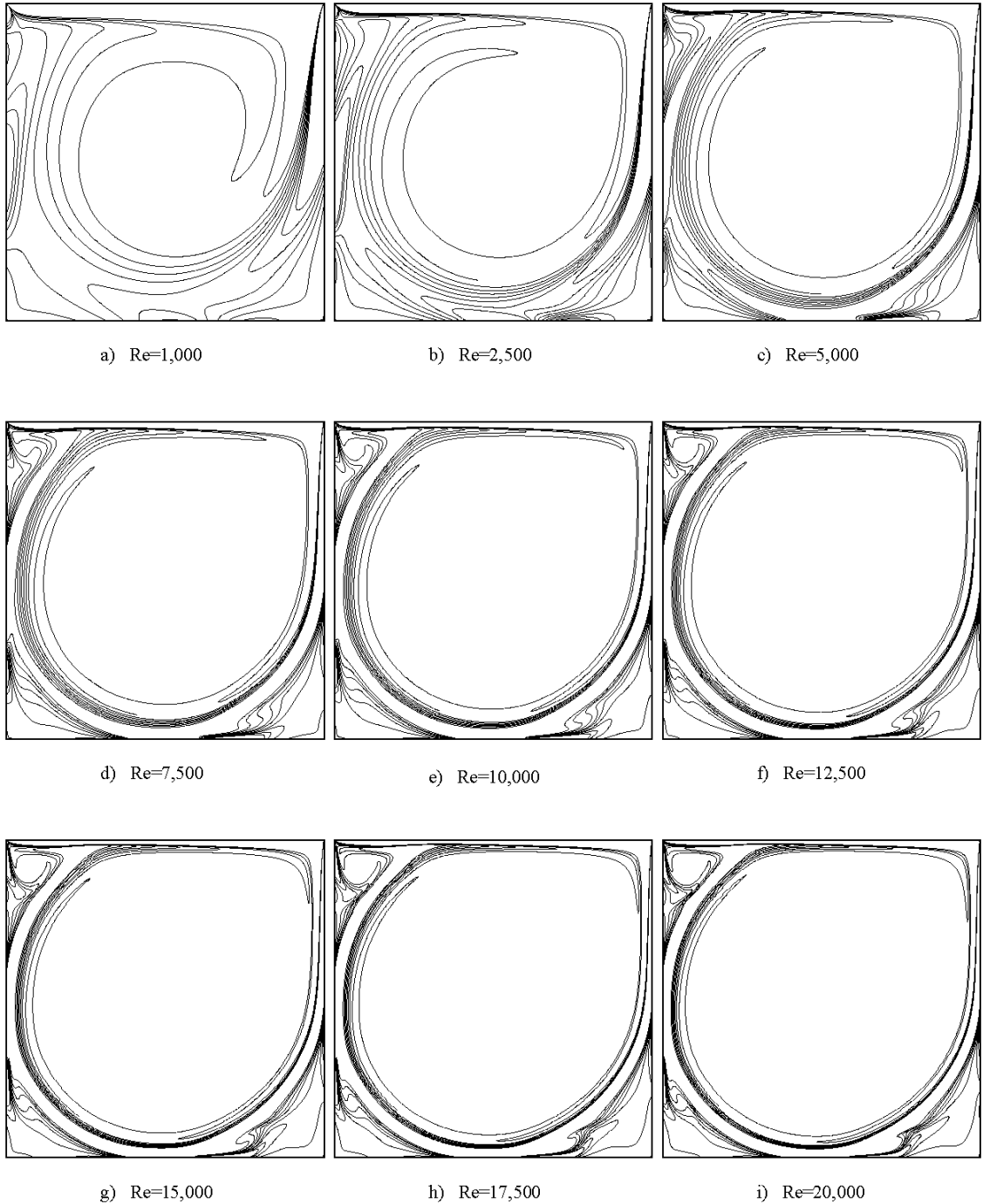

Figure 2. Vorticity contours at various Reynolds numbers: (a) Re=1000;
(b) Re=2500; (c) Re=5000; (d) Re=7500; (e) Re=10000; (f) Re=12500;
(g) Re=15000; (h) Re=17500; and (i) Re=20000.

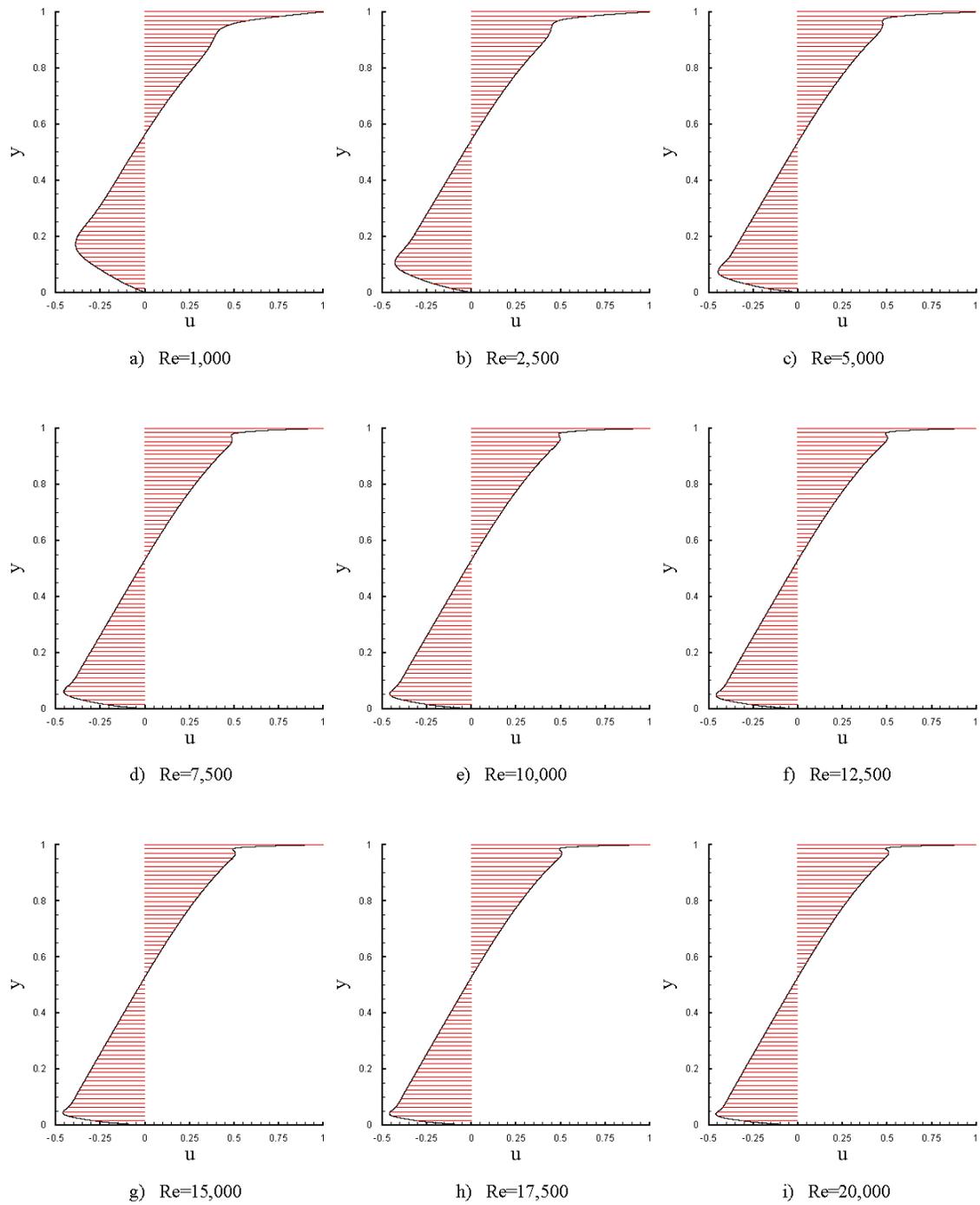

Figure 4. The u-velocity profiles along a vertical line passing through the centre of the cavity:
(a) Re=1000; (b) Re=2500; (c) Re=5000; (d) Re=7500; (e) Re=10000;
(f) Re=12500; (g) Re=15000; (h) Re=17500; and (i) Re=20000.

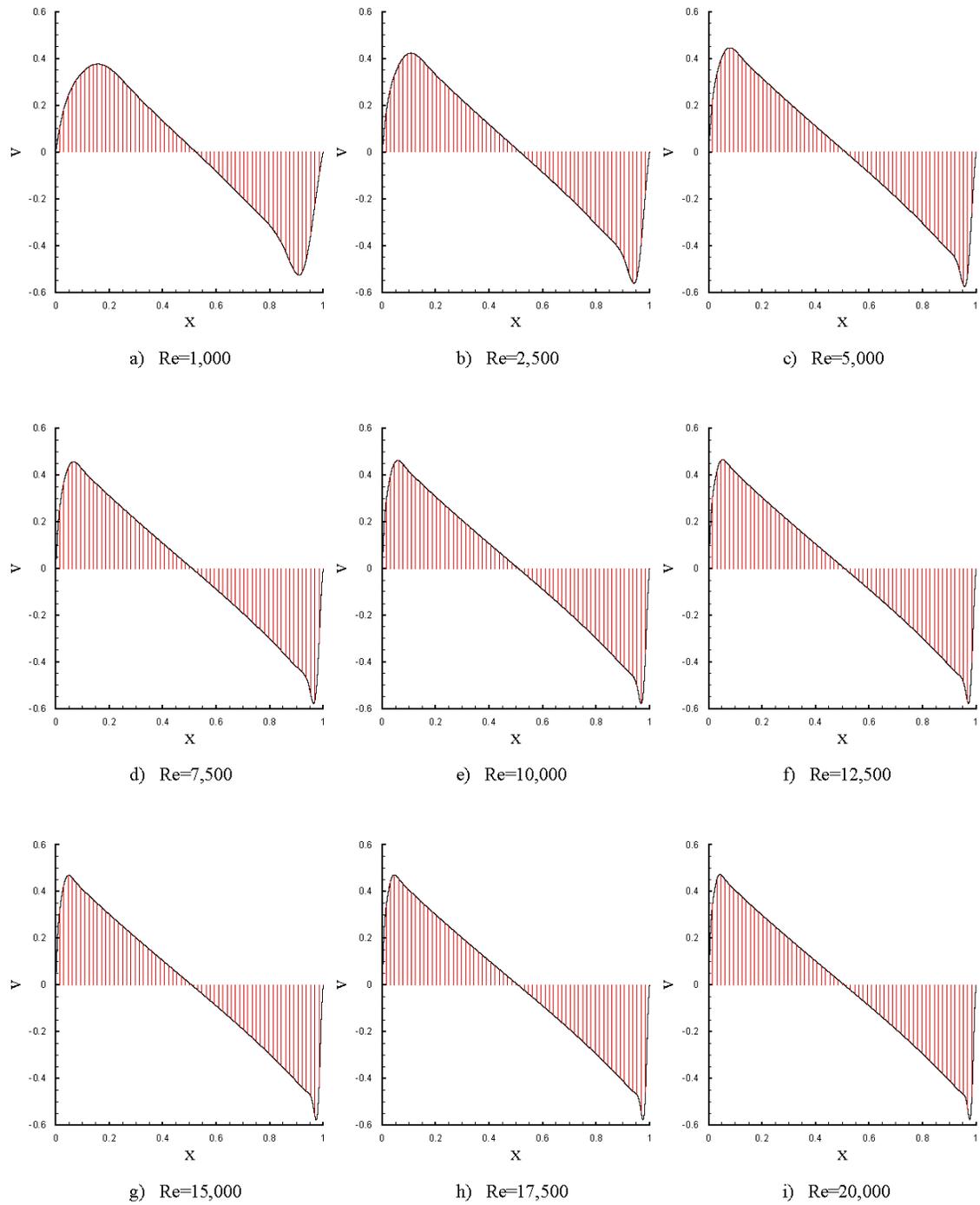

Figure 5. The v-velocity profiles along a horizontal line passing through the centre of the cavity:
(a) Re=1000; (b) Re=2500; (c) Re=5000; (d) Re=7500; (e) Re=10000; (f) Re=12500;
(g) Re=15000; (h) Re=17500; and (i) Re=20000.

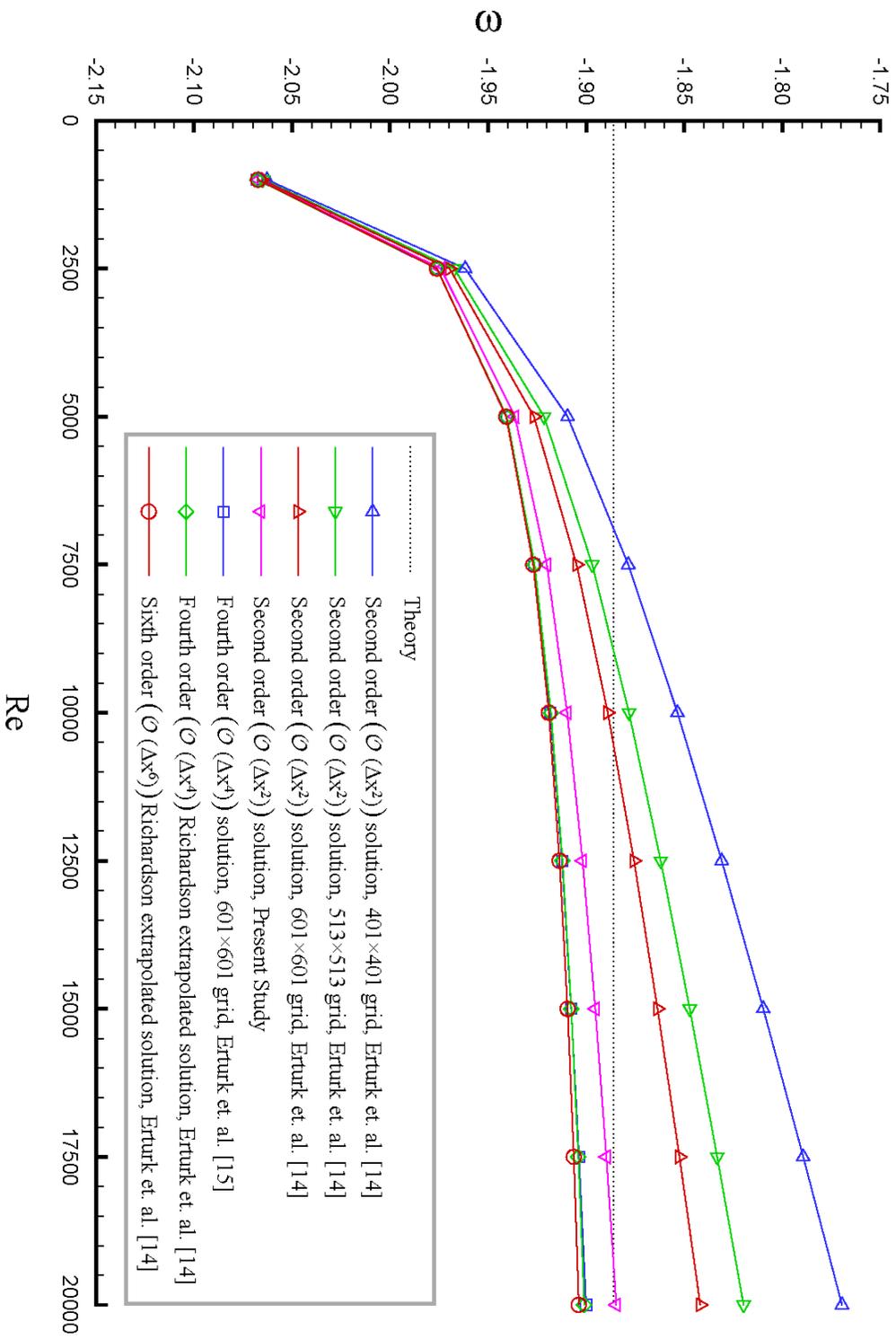

Figure 6. Comparison of numerical solutions with theoretical solution.